\documentclass[runningheads,a4paper,12pt]{llncs}
\usepackage{amssymb}
\usepackage{graphicx}
\usepackage{amsmath}
\usepackage{microtype}

\newcommand{\R}{\mbox{$\mathbb{R}$}}

\newcommand{\C}{\mbox{$\mathbb{C}$}}
\newcommand{\N}{\mbox{$\mathbb{N}$}}

\newcommand{\Z}{\mbox{$\mathbb{Z}$}}

\newcommand{\Q}{\mbox{$\mathbb{Q}$}}
\newcommand{\F}{\mbox{$\mathbb{F}$}}
\newcommand{\A}{\mbox{$\mathbb{A}$}}

\newcommand{\keywords}[1]{\par\addvspace\baselineskip
\noindent\keywordname\enspace\ignorespaces#1}

\begin{document}

\mainmatter  

\title{Hensel's Lemma, Backward Dynamics and $p$-adic Approximations}

\titlerunning{Backward Dynamics and $p$-adic Approximations}

%
%
\author{  Sushma Palimar}

%
\authorrunning{Backward Dynamics and $p$-adic Approximations }

\institute{Department of Mathematical and Computational Sciences\\
National Institute of Technology Karnataka, Surathkal, India.\\
}
%
%

\toctitle{Backward Dynamics}
\tocauthor{$p$-adic Approximations}
\maketitle

\begin{abstract}
The problem of backward dynamics over the ring of $p$-adic integers is studied. It is shown that 
Inverse Limit Theory provides the right framework. Backward iterations of a polynomial with $p$-adic integer coefficients
are constructed by solving congruences modulo powers of $p$, which inturn are solved by Hensel's lifting lemma.

\emph
\end{abstract}

\section{Introduction}\label{sect-1-9-12-11}
Dynamical systems originally arose in the study of systems of 
differential equations used to model physical phenomena.
One simplification in this study is to discretize time, so
that the state of the system is observed only at discrete steps of time.
 This leads to the study of the iterates of a single
transformation. One is interested in both quantitative behavior, such as
the average time spent in a certain region, and also qualitative behavior,
such as whether a state eventually becomes periodic or tends to infinity.\\
A discrete-time dynamical system consists of a non-empty set $X$ and a map
$f:X\rightarrow X$. 
For $n\in \N$,\label{subpar3.1} the $n$th iterate of $\phi$ is the $n$-fold 
composition $\phi^{[n]}=\phi\circ \phi\circ...\circ \phi\,$ ($n$ times); $\phi^{[0]}$
is defined to be the identity map. If  $\phi$ is invertible, 
then $\phi^{-n}=\phi^{-1}\circ \phi^{-1}\circ...\circ \phi^{-1}$ 
($n$ times). Since $\phi^{[n+m]}=\phi^{[n]}\circ \phi^{[m]}$, these iterates 
form a group if $\phi$ is invertible, and a semigroup otherwise.
For a given $\alpha$ in $X,$ the forward orbit of $\alpha$ is the set
\[{\cal O_{\phi({\alpha})}}=\{\phi^{n}({\alpha}):n\geq0\}\]
If the orbit ${\cal O_{\phi({\alpha})}}$ is finite then $\alpha$ is said to
be a  pre-periodic point, otherwise $\alpha$ is said to be a wandering point.\\
The central problem in dynamics is to classify the points $\alpha$ in the set $X$
according to the behavior of their orbits  ${\cal O_{\phi({\alpha})}}$.
In practice $X$ usually has additional structure that is preserved 
by the map $f$. For example, $(X,f)$ could be a measure
space and a measure preserving map; a topological space and a continuous map; 
a metric space and an isometry; or  a smooth manifold and a differentiable map;
a finite set (e.g., Finite Field) and a polynomial.
\subsection{Symbolic dynamics} Symbolic dynamics arose as an attempt to study such systems by means
of discretizing space as well as time. The basic idea is to divide up the set
of possible states into a finite number of pieces. Each piece
is associated with a ``symbol'', and in this way the evolution of the system
is described by an infinite sequence of symbols. This leads to a ``symbolic"
dynamical system that mirrors and helps us to understand the dynamical
behavior of the original system. Computer simulations of continuous systems 
necessarily involve a discretization of space, and results of symbolic dynamics help us
understand how well, or how badly, the simulation may mimic the original.
 Symbolic dynamics by itself has proved a bottomless source of beautiful mathematics
and intriguing questions. As polygons and curves are to geometry shift spaces
are to symbolic dynamics. The set
\[\Sigma=\{0,1,2,...,m-1 \}^{\N}\]
is called the sequence space on $m$ symbols $0,1,...,m-1$.\\
The most important ingredient in the sequence space is the 
shift map $\sigma.$
The shift map $\sigma: \Sigma\rightarrow\Sigma$
is given by $\sigma((s_0,s_1,s_2,..))\,=\,(s_1,s_2,s_3,...)$.\\
The shift map discards the first entry in the sequence 
and shifts all other entries one place to the left. \\
The distance between two sequences $s=(s_0,s_1,...)$ and $t=(t_0,t_1,...)$ is 
given by $d(s,t)=\sum_{i=0}^{\infty}{\frac{s_i-t_i}{m^{i}}}$. \\
For an integer $m>1,$ set $\A_{m}=\{1,2,...,m\}$.
Let $\Sigma_{m}={\A^{Z}_{m}}$ be the set of infinite two sided 
sequences of symbols in  
$\A_{m}$ and $\Sigma^{+}_{m}={\A^{N}_{m}}$ be the set of infinite one-sided sequences.
The pair $(\Sigma_{m},\sigma)$ is called the full two sided shift; 
$(\Sigma^{+}_{m},\sigma)$ is called the full one sided shift.
The two-sided shift is invertible. For a one-sided sequence, the leftmost symbol disappears,
so the one-sided shift is non-invertible, and every point has $m$ pre-images.
Both shifts have $m^{n}$ periodic points of period $n$.
The shift spaces are compact topological spaces in the product topology.
This topology has a basis consisting of cylinders
$C^{n_{1},...,n_{k}}_{j_{1},..,j_{k}}$=$\{x=(x_l):x_{n_{i}}=j_{i}, i= 1,2,...,k\}$, 
$n_1<n_2<...<n_k$ are indices in $\Z$ or  $\N$, and $j_{i}\in \A_{m}$.
 Since the preimage of a cylinder is a cylinder, $\sigma$ is continuous on
 $\Sigma^{+}_{m}$ and homeomorphism on $\Sigma_{m}$. The metric 
$d(x,x^{\prime})=2^{-l}$, for $l=min\{|i|:x_i\neq x_{i}^{\prime}\} $
generates the product topology on $\Sigma^{+}_{m}$ and $\Sigma_{m}$.
In the product topology periodic points are dense and hence there are dense orbits.
A detailed study can be found in \cite{Dev}, \cite{Brin}. 
\section{Arithmetic Dynamical Systems}
Classically, discrete dynamics refers to 
the study of the iteration of self-maps of the complex plane or real line.
Arithmetic Dynamics is the study of number theoretic properties of dynamical systems.
Arithmetic dynamics is discrete-time dynamics (function iteration) over arithmetical sets, such as
algebraic number rings and fields, finite fields, 
$p$-adic fields, polynomial rings, algebraic curves,
etc. A thorough introduction  is given in  \cite{silver}, \cite{silver1}.\\
 In this paper we study the problem of backward 
dynamics of any polynomial of finite degree over finite rings 
$\Z/p^{n}\Z$ using Hensel's lifting lemma, as usual $p$ denotes a prime.\\
Backward-iteration sequences  given by\[x_{n}=f(x_{n+1}),\qquad n>0\] are of a different 
nature because a point could have infinitely many pre-images as well as none.
If the given forward moving map is a quadratic map, 
the corresponding backward map is a square root map; if the given map is a cubic map the corresponding 
backward map is the cube root map and so on. Thus 
essentially we solve  $f(x)=0$ as a polynomial over the defining set.
For e.g., the Julia set can be found as the set 
of limit points of the set of pre-images of (essentially) any given point.
Unfortunately, as the number of iterated pre-images grows
exponentially, this is not feasible computationally
when the underlying set is the set of Real or Complex numbers.\\
In general, maps of higher degree ($\geq 5$) 
are not suitable  for backward dynamics over $\R$ or $\C$.
As there is no explicit formula for solving  a polynomial of degree $\geq 5$,
roots can be found by using standard techniques from  Nummerical Methods    over $\R$ or $\C$
and we either arrive at a null sequence or constant sequence after some backward iteartion.
In such cases nothing can be said about the behavior of trajectories. 
But when the set is finite (endowed with algebraic structure) 
and  the map is a polynomial it is possible to retrieve 
the pre-images (roots), if they exist, with respect to different 
prime power moduli and thus study the structure of pre-images locally at that prime.\\
We show that this problem can be well understood over the ring of $p$-adic integers $\Z_{p}$.
This process deals with an imporatnt branch of mathematics called \textit{Inverse Limit Theory}.
Below we discuss some basics of \textit{Inverse Limit Theory} and the $p$-adic integers.
\subsection{Inverse Limits}
Let $X_0,X_1,X_2,...$ be a countable collection of spaces, and suppose that, for each $n>0,$
there is a continuous mapping $f_n:X_n\rightarrow X_{n-1}$. 
The seqence of spaces and mappings $\{X_n,f_n\}$ is called an \textit{inverse limit sequence}
and may be represented as 
\[ ...\xrightarrow {f_{n+1}}X_{n}\xrightarrow{f_{n}} X_{n-1}...\xrightarrow{f_{2}} X_{1}\xrightarrow{f_{1}}{X_0} \]
Clearly, if $n>m,$ there is a continuous mapping $f_{n,m}:X_{n}\rightarrow X_{m}$
given by the composition
$f_{n,m}=f_{m+1}\cdot f_{m+2}\cdot \cdot \cdot f_{n-1}\cdot f_{n}$.\\
Consider a sequence $(x_0,x_1,...,x_n,..)$ such that each $x_{n}$ is a point of the space  $X_{n}$ 
and such that  $x_n= f_{n+1}(x_{n+1}), \; n\geq 0.$
Such a sequence can be identified in the product space 
$\prod_{n=0}^{\infty} X_{n}$ by considering a function $\phi$
from the nonnegative integers into $\prod_{n=0}^{\infty}X_{n}$,
given by $\phi(n)=x_{n}$. Thus the set of all sequences is a subset of 
$\prod_{n=0}^{\infty} X_{n}$ and has a topology as a subspace. This 
topological space is the inverse limit space of of the sequence $\{ X_n,f_n\}$
denoted by $X=\varprojlim (X_{n},f_{n})$.  
\begin{theorem} \label {8-12-th.2}
 Suppose that each space $X_{n}$ in the inverse limit sequence $\{ X_{n},f_{n}\}$ is a
compact Hausdorff space. Then $X$ is not empty \cite{gill}.
\end{theorem}
\begin{theorem} 
 A space $X$ is a compact Hausdorff space with $dim(X)\leq 0$ if and only if $X$
is an inverse limit of finite discrete spaces \cite{nagem}. 
\end{theorem}
A finite discrete space is totally disconnected, 
compact and Hausdorff and all those properties carry over to inverse limits too.\\
A detailed study of Inverse limit spaces can be found 
in \cite{gill}.
\subsection{$p$-adic Integers}
An important class of such inverse limits is given by rings of $p$-adic integers $\Z_{p}$.
For every $n\geq 1$, let $A_n= \Z/p^{n}\Z$. An element of $A_{n}$ defines in an 
obvious way an element of $A_{n-1}$ and the homomorphism 
\[\phi_{n}:A_{n}\rightarrow A_{n-1}\] is surjective and the kernel is $p^{n-1}A_{n}$.  
The sequence 
\[ ...\rightarrow A_{n}\rightarrow A_{n-1}...\rightarrow A_{2}\rightarrow{A_1} \]
forms a ''projective system''  indexed by the integers $\geq 1$. 
Inverse limit of this inverse system is  $\Z_{p}=\varprojlim (A_{n},\phi_{n})$. Refer  \cite{serre} for details.
\begin{definition}
 The ring of $p$-adic integers, $\Z_{p}$,  is the projective limit (inverse limit) of the system $(A_{n},\phi_{n})$.
\end{definition}
An element of $\Z_{p}=\varprojlim (A_{n},\phi_{n})$ is a sequence $x=(...,x_n,...,x_1)$
with  $x_n\in A_n$ and $\phi_n({x_n})=x_{n-1}$ if $n\geq 2$.
Addition and multiplication in $\Z_{p}$ are defined co-ordinate wise.
In other words, $\Z_p$ is a subring of the product $\prod_{n\geq 1}A_{n}$.
If $A_{n}$ is endowed with discrete topology and $\prod_{n\geq 1} A_n$ the product topology,
the ring $\Z_p$ inherits a topology which turns it into a compact space. \\
Let $p$ be a prime. For $n\in \Z$,  let $\nu_{p}(n)$ denote the exponent of highest power of $p$ that 
divides $n$, and  $\nu_{p}(0)\,=\infty$. More formally, $\nu_{p}(n)$   is the 
unique natural number such that $n\,=\,p^{\nu_{p}(n)}u$ with $p\nmid u.$ This definition can be extended 
to $\Q$ by letting $\nu_{p}(\frac{a}{b})\,=\nu_{p}(a)-\nu_{p}(b)$. The $p$-adic
 absolute value is defined by $|x|_{p}=p^{-\nu_{p}(x)}$ for any $x\in \Q$.
 Then $d_{p}(x,y) = |x-y|_{p}$ defines a metric  on $\Q$.
The metric space $(\Q,d_{p})$ is not complete, and its completion is the $p$-adic number field $\Q_{p}$.
This absolute value is \textit{non-Archimedean} as, in the place of triangle inequality,
the stronger relation $|x+y|_{p}\leq\,Max\{ |x|_{p},|y|_{p}\}$,
 also known as \textit{ultrametric inequality} holds. 
The non-Archimedean property is equivalent to the assertion that, 
 sup$\{|n|\,: n \in \Z\}=1$\\
As a consequence (in stark contrast to the $Euclidean$ norm), 
the $p$-adic norm does not permit accumulation of errors in the following sense: if each of $k$ elements 
$\{x_{1},x_{2},...,x_{k}\}$ have $p$-adic norm atmost $\epsilon$,
 then $|x_{1}+x_{2}+...+x_{k}|_{p}\leq\epsilon$ as well.
This property justifies extensive use of modular 
arithmetic (\textit{$p$-adic estimation}) in $p$-adic calculations.
\begin{definition}\label{def-9-12-11}
 A $p$-adic integer is a formal series $\sum_{i\geq 0}a_{i}p^{i}$ with integral coefficients $a_i$
satisfying $0\leq a_{i}\leq p-1.$
\end{definition} 
The subset of $\Q_{p}$ defined by $\Z_{p}=\{ x\in \Q_{p}:|x|_{p}\leq1\}$ 
is called the \textit{integer ring} and is an integral domain. 
The set of invertible elements in $\Z_{p}$, called its group of units 
is $\Z_{p}^{*}=\{x\in\Z_{p}:|x|_{p}=1 \}$. The ring $\Z_{p}$ contains a unique
maximal ideal $p\Z_{p}=\{x\in\Z_{p}:|x|_{p}<1 \}$.
 The quotient of $\Z_{p}$ by this maximal ideal is a field,
 which can be identified with the finite field of $p$ elements
 $\F_{p}=\Z/p\Z$ in the obvious way. $\F_{p}$ is called the
 \textit{residue class field} of $\Q_{p}$. The field $\Q_{p}$ is
unorderable, in essence due to the modular arithmetic of $\F_{p}$. 
It has characteristic $0$ since it contains $\Q$ as
 a subfield. Indeed, $\Q$ is a dense proper subset of $\Q_{p}$ and embeds into
 $\Q_{p}$ as the set of elements whose $p$-adic 
coefficients are eventually periodic. 
Topologically $\Q_{p}$ is a Cantor set: totally disconnected
 but not discrete. In algebraic sense too $p$-adics 
are full of holes: There is no finite extension of $\Q_{p}$ which is algebraically closed. 
Ostrowski's theorem states that any non-trivial 
absolute value on the rational numbers $\Q$ is equivalent 
to either the usual real absolute value or a $p$-adic absolute value.
 Refer \cite{Borevich}, \cite{Robert} for details.
\section{$p$-adic Dynamics }
It is natural to have accumulated truncation errors or round-off errors
in a  dynamical system even with a small perturbation. 
These are unavoidable since even with simple repetitive operations,
the number of digits of the result can increase so much that 
the result cannot be held fully in the registers 
available in the computer. Such errors accumulate one after 
another from iteration to iteration generating new errors. 
These difficulties motivated to look for an alternate number 
system which possesses the best features as well as the 
advantages of both the $p$-ary and residue number system. 
Such a number system is the $p$-adic number system, 
discovered by Kurt Hensel in 1897 in the course of his 
work on finding new completions of the rational numbers. 
Hensel's original description of the $p$-adic numbers involved an analogy between
 the ring of integers and the ring of polynomials over the complex
numbers, the crux of which was the development of a representation of rational 
numbers analogous to that of Laurent expansions of rational functions
namely, the $p$-adic  expansion.
This idea was motivated by the 
existence of real expansions of rational numbers with respect to a $p$-scale:
\begin{displaymath}
 x=\sum_{n=-\infty}^{k} \alpha_{n}p^{n}, \qquad \alpha_{n}=0,...,p-1.
\end{displaymath}
Such manipulations with rational numbers and series generated the idea that
there exists some algebraic structure similar to the system of real numbers 
$\R$. Thus each $\Q_{p}$
has the structure of a number field. In fact, the fields of $p$-adic numbers,
$\Q_{p}$, were the first examples of infinite fields that differs from 
$\Q$, $\R$, $\C$ and corresponding fields of rational functions.
The following definition is useful in identifying a $p$-adic integer.
\begin{definition} 
 Let $p$ be some prime number. A sequence of integers 
\[\{x_n\}=\{x_0,x_1,...,x_n,...\}\]
satisfying
\[x_n\equiv x_{n-1} \pmod{p^{n}}\]
for all $n\geq 1,$ is called a $p$-adic integer.
Two sequences and $\{x_n\}$ and $\{x_{n}^{\prime}\}$ determine the same $p$-adic integer
if and only if $x_{n}\equiv x_{n}^{\prime}\pmod{p^{n}}$ for all $n\geq 0$
\end{definition}
This definition can be easily identified with the definition \ref{def-9-12-11}.
\paragraph*{}
The connection between congruences and equations is based on the simple 
remark that, if the equation \[ F(x_1,x_2,...,x_n)=0\]
where $F$ is a polynomial with integral coefficients, has a solution in integers,
then the congruence \[F(x_1,x_2,...,x_n)\equiv 0 \pmod{m}\]  
is solvable for any value of the modulus $m$. 
On the other hand the situation is more complicated for congruences.
For any modulus $m>1$, there are polynomial congruences having no solutions.
For eg., the congruence $x^{p}-x+1\equiv 0 \pmod{m}$ has no solution if 
$p$ is any prime factor of $m$ by Fermat's theorem. Where as a congruence
can have more solutions  than its degree, for eg., $x^{2}-7x+2\equiv 0\pmod{10}$
has four solutions $x=3,4,8,9.$ But if the modulus is a prime, a congruence
cannot have more solutions than its degree. \\
The following results are basic, see \cite{Zucker}. 
 \begin{theorem}\label{7-12.th1}
If the degree $n$ of $f(x)\equiv 0 \pmod{p}$ is greater than or equal to $p$,
then either every integer is a solution of $f(x)\equiv 0 \pmod {p}$ or there 
is a polynomial $g(x)$ having integral coefficients, with leading coefficient $1$,
such that, $g(x)\equiv 0 \pmod{p}$ is of degree less than $p$ and the solutions 
of  $g(x)\equiv 0 \pmod{p}$ are precisely those of $f(x)\equiv 0 \pmod{p}$.
\end{theorem}
\begin{theorem}\label{7-12.th2}
 The congruence $f(x)\equiv 0 \pmod{p}$ of degree $n$ has atmost $n$ solutions.
\end{theorem}
\begin{corollary}\label{7-12.co1}
 If $b_{n}x^{n}+b_{n-1}x^{n-1}+...+b_{0}\equiv 0 \pmod{p}$ has more than $n$ solutions,
then all the coefficients $b_{j}$ are divisible by $p$.
\end{corollary}
\begin{theorem}\label{7-12.th3}
 The congruence $f(x)\equiv 0 \pmod{p}$ of degree $n,$ with leading coefficient
$a_{n}=1,$ has $n$ solutions if and only if $f(x)$ is a factor of $x^{p}-x$ modulo $p,$
that is, if and only if $x^{p}-x=f(x)q(x)+ps(x)$, where $q(x)$ and  $s(x)$ have integral
coefficients, $q(x)$ has degree $p-n$ and leading coefficient $1,$ and where $s(x)$ is a polynomial
of degree less than $n$ or $s(x)$ is zero. 
\end{theorem}
The proofs of the above theorems are  simple consequences of Fermat's little theorem and its application 
which can be found in many books on number theory. 
\subsection{Hensel's Lifting Lemma}
As both $\Q$ and $\R$ are  not algebraically closed, they do not always contain all 
roots of polynomials with integer coefficients.
 Though $\C$ is algebraically closed, as the degree of the polynomial increases, finding roots of the polynomial 
is computationally not feasible, as it requires high level of precision for machine computation.
To this end, we now turn our attention to solving polynomial congruences
modulo prime powers.\\ We note that for any polynomial $f(x) \in  \Z[x]$ and any
integer $r$, there is a polynomial $g_{r}(x) \in  \Z[x]$ with
$f(x + r) = f (r) + x f^{\prime}(r)+x^{2}g_{r}(x)$.\\
This can be seen either through the Taylor expansion for $f(x + r)$ or through
the binomial theorem in the form\begin{equation} \label{eq7-12.1}
                                (x + r)^{d}=r^{d} + dr^{d-1}x+x^{2}\sum_{j=2}^{d}\binom{d}{j}r^{d-j}x^{j-2} 
                                \end{equation}
The above standard results mentioned in the form of Theorems can be used to find the  solutions to $f(x)\equiv 0 \pmod{p}$.
The question is how we might be able to “lift” a solution to one modulo $p^{k}$
 for various exponents $k$. With regard to this, 
Hensel's Lemma is a powerful tool which relates the roots of a given polynomial to its solution modulo a
prime. The lemma and its proof both rely on iterative procedures that return an agreeable solution
if supplied with a well-behaved seed. 
\begin{definition}
 If $f(a)\equiv 0 \pmod{p}$, then the root $a$ is called $nonsingular$  if $f^{\prime}(a) \not\equiv 0\pmod{p}$;
otherwise it is singular.
\end{definition}
Two versions of \textit{Hensel's Lemma} are stated below.
\begin{theorem}\label{th-11-12-11.1}
 Hensel's Lemma over the ring of integers.\\
  Suppose that $f(x)$ is a polynomial with integral coefficients.
 If $f(a)\equiv 0\pmod{p^{j}}$ and 
$f^{\prime}(a)\not\equiv 0 \pmod{p}$, then there is 
a unique $t\pmod{p}$ such that \\$f(a+tp^{j})\equiv 0\pmod {p^{j+1}}$.
\end{theorem}
\begin{theorem}\label{th-11-12-11.2}
 Hensel's Lemma over the ring of $p$-adic integers.\\
 Let $f\in   \Z_{p} [x]$ be monic. If $a_{0} \in  \Z$ is a simple root of
$f (x)\equiv 0 \pmod{p}$,
then $\exists\, y \in \Z_{p}$ such that $y\equiv  a_{0} \pmod {p}$ and $f(y) = 0.$
\end{theorem}
\begin{proof}
Suppose that $\exists a_{n}$ such that $f(a_{n} ) \equiv  0 \pmod{p^{n}}$.
We must show that $a_{n}$ can be lifted uniquely to $a_{n+1} \pmod{p^{n+1}}$ such that
 $a_{n+1}\equiv a_{n} \pmod{p^{n}}$ and $f(a_{n+1} )\equiv 0 \pmod {p^{n+1} }$,
 then $y$ is the limit of this sequence of
$\pmod {p^{k}}$ solutions.\\
Since $f$ is a polynomial we can write it in the form $f(x)= \sum_{i}c_{i}x^{i}$.
Also consider $tp^{n}+a_{n}$
as a possible lift of $a_{n}$ . Then
\begin{align}
 f(a_n+tp^{n})=\sum_{i}c_{i}(tp^{n}+a_{n})^{i}\\
\equiv f(a_{n} )+p^{n}tf^{\prime}(a_n)\pmod {p^{n+1}}.
\end{align}
The equivalence above is a result of Taylor series expansion.
Now, solve for $t$ in 
\[
p^{n} tf^{\prime}(a_{n})+ f(a_{n})\equiv 0\pmod {p^{n+1}}.\]
Thus \[tf^{\prime}(a_{n} )\equiv -\left( \frac{f (a_n)}{p^{n}}\right) \pmod{p}.\] 
Since $f(a_{n} ) \equiv 0 \pmod{p^{n}}$, (since $a_{n} \equiv a_0 \pmod {p})$
and $f^{\prime}(a_n)\not\equiv 0 \pmod {p}$ (simple root), then $t$ has a unique solution
$\pmod {p}$. Thus
$a_{n+1} = a_{n} + tp^{n}$ is a unique lift of $a_{n} \pmod {p}$.
Thus we have constructed an infinite sequence of $a_{i}$ such that $f(a_{i}) \equiv 0 \pmod{p^{i}}$, such that, 
$f^{\prime}(a_{i}) \not\equiv0 \pmod {p^{i}}$  
and $a_{i+1} \equiv a_{i} \pmod {p}$. This sequence is Cauchy, and therefore converges to a
unique limit $y \in  \Z_{p}$.
\end{proof}
Since $\Z$ is dense in $\Z_{p}$, proof of Theorem \ref{th-11-12-11.1} follows directly from  Theorem \ref{th-11-12-11.2}.\\ 
It can be verified that, proof of 
Hensel's Lemma is entirely analogous to Newton's method for locating 
the root of a differentiable function.
Let us recall Newton's method from calculus as a method  of finding roots to a polynomial by choosing
a seed and then making better and better approximations based on the polynomial's derivative at
that point. In the case of Newton's method, the condition on the seed is that the derivative at that
point be non-zero, otherwise it supplies no useful information for improving at each iteration.
   Hensel's Lemma is similar, it takes a polynomial with coefficients in $\Z_{p}$ and instead of
requiring a ‘guess’ at a possible root, it requires a $p$-adic integer that is a root mod ${p}$, 
i.e. some $\alpha$ such that the polynomial $f(x)$ evaluated at $\alpha$ is \[f(\alpha)\equiv 0 \pmod{p\Z_{p}}\]
This method will then return roots mod $p, p^{2} , p^{3} , . . . $
until the desired root of the equation is found. 
\subsection{Backward Iterations, Inverse Limits and $p$-adic Approximations }
Now we find the backward iteration of any given polynomial $f(x)$ of arbitrary degree say $n$,
at any point, say $x_n$ of the forward iterating orbit, 
i.e., we would like to solve \begin{equation}\label{eq8-12.1} f(x)\equiv x_n\pmod{p^{k}}\end{equation} for some $k\in \N$.
For this, we first solve the congruence $f(x)\equiv x_{n} \pmod{p}$,  such that the nonzero coefficients of $f(x)$
are relatively prime to $p$ and $x_{n}$ is chosen as above,  so that the corresponding 
backward iterating orbit of $x_{n}$ are found modulo $p^{j}$, $j\rightarrow \infty.$ \\
If the degree $n$ of $f(x)$ is greater than $p,$ then by Theorem \ref{7-12.th1}, $f(x)$ is divided by
($x^{p}-x$) mod $p$ and   solutions of the resulting polynomial $g(x)$ are  the same as  
those of $f(x)\equiv x_{n} \pmod{p}$. Since the modulus is prime, the 
congruence cannot have more solutions than its degree. \\
Let $a_1,a_2,...,a_l$ be the roots obtained, by solving the  congruence $f(x)\equiv x_{n} \pmod{p}$
 where 
($l\leq deg(g(x)$ by Theorem \ref{7-12.th2}). Now each $a_{i}$ is lifted modulo $p^{j}$, $j=2,3,...$
whenever $f^{\prime}(a_{i})\not\equiv 0 \pmod{p}$, ($i=1,2,..,l.$)i.e., whenever the $a_{i}'s$ are nonsingular.
We work up to a fixed precision say $j=k$.
Once the roots are lifted modulo $p^{k}$, $x_n$ is replaced by the lifted root and
the congruence equation (\ref{eq8-12.1})
is soved now with respect to the new root.
The process of replacing the old root by the \textit{new lifted root} is repeated for some finite number of steps.
Thus with respect to each nonsingular $a_{i}$, we obtain the corresponding backward sequence, generated
from the single seed $x_{n}$.  Thus if there are $n$ nonsingular roots, after 
 the $m$th step of replacing the old root by the \textit{new lifted root},
there are atmost $n^{m}$ backward iterating points generated from the single seed $x_{n}$.
Hence there exists a tree like structure, the roots may be called as leaves and the  branches are formed
at each new lifted root.
A simple code for this program can be written on python, which computes the sequences
effectively.\\  Sequences thus obtained belong to the \textit{inverse limit space} by definition.
It can be verified that the sequence space formed by the above backward iterations 
 is totally disconnected and discrete.  
Also, if the original spaces are discrete, then the inverse limit space 
  $\varprojlim \{X,\,f\}$  is totally disconnected. 
This is one way of 
realizing the $p$-adic numbers and the Cantor set. \\To study the long time behaviour of a dynamical system 
it is necessary to introduce a suitable metric. 
A natural choice would be the one given in section \ref{sect-1-9-12-11}, i.e., 
the distance between two sequences $s=(s_0,s_1,...)$ and $t=(t_0,t_1,...)$ is 
given by $d(s,t)=\sum_{i=0}^{\infty}{\frac{s_i-t_i}{p^{i}}}$. \\
Thus by the introduction of a metric, sequence space of backward iteration becomes a compact metric space.
The sequence space thus obtained can be identified with the ring of 
$p$-adic integers in view of the following theorems:
\begin{theorem}\label{thgil}
 Any two totally disconnected perfect compact metric spaces are homeomorphic.
\end{theorem}
\begin{theorem}\label{8-12-th.1}
 Let $M$ be a compact , totally disconnected metric space. Then $M$ is homeomorphic to 
the inverse limit space of an inverse limit sequence of finite discrete spaces \cite{gill}.
\end{theorem}
\section{Conclusion}
The problem of backward iteration is studied by solving congruences, 
which in turn are solved by \textit{Hensel's lifting}.
Sequences generated by such solutions form, naturally, elements of an  Inverse Limit Space. 
These spaces have  also been characterised.

\end{document}